\documentclass[12pt]{amsart}
\usepackage{amssymb}
\usepackage{colordvi}
\usepackage{graphicx}
\usepackage{vmargin}
\setpapersize{USletter}
\setmargrb{1in}{1in}{1in}{1in} % --- sets all four margins.  Cool.
\newtheorem{theorem}{Theorem}[section]
\newtheorem{proposition}[theorem]{Proposition}
\newtheorem{lemma}[theorem]{Lemma}

\newtheorem{defn}[theorem]{Definition}

\theoremstyle{definition}

\newtheorem{example}[theorem]{Example}
\newtheorem{remark}[theorem]{Remark}
\newenvironment{ex}{\begin{example}\rm}{\end{example}}

\newcommand{\R}{\mathbb{R}}

\newcommand{\PP}{\mathbb{P}}
\newcommand{\TP}{{\mathbb T}{\mathbb{P}}}

\title{Tropical quadrics through three points}

\author{Sarah~B.~Brodsky}
\address{Sarah~B.~Brodsky \\
Universit\"at Kaiserslautern \\
 Fachbereich Mathematik \\
 Erwin-Schr\"odinger-Stra\ss e \\
 D-67633 Kaiserslautern, Germany}
\email{brodsky@mathematik.uni-kl.de}

\author{Bernd Sturmfels}
\address{Bernd Sturmfels\\
Department of Mathematics\\
      University of California\\
      Berkeley, California 94720,
      USA}
\email{bernd@math.berkeley.edu}
\begin{document}

\begin{abstract}
We tropicalize the rational map that takes triples of points in the 
projective plane to the plane of quadrics passing through these points.
The  image of its tropicalization is contained in the tropicalization of 
its image. We identify these objects inside the
tropical Grassmannian of planes in projective $5$-space,
and we explore a  small tropical Hilbert~scheme. 
\end{abstract}

\maketitle

\section{Introduction}
Given three points $x = (x_0{:}x_1{:}x_2)$,
$y = (y_0{:}y_1{:}y_2)$ and $z = (z_0{:}z_1{:}z_2)$
in the projective plane $\PP^2$ over a field $K$, we are interested in
the space $L_{x,y,z}$ of all homogeneous quadrics
that vanish at $x$, $y$ and $z$. By definition, the vector space
$L_{x,y,z}$ is the kernel of the $3 \times 6$ matrix
\begin{equation}
\label{3by6matrix}
\begin{pmatrix}
x_0^2 & x_1^2 & x_2^2 & x_0 x_1 & x_0 x_2 & x_1 x_2 \\
y_0^2 & y_1^2 & y_2^2 & y_0 y_1 & y_0 y_2 & y_1 y_2 \\
z_0^2 & z_1^2 & z_2^2 & z_0 z_1 & z_0 z_2 & z_1 z_2
\end{pmatrix}.
\end{equation}
This defines a rational map which is a morphism on the open set of
non-collinear triples:
$$\Phi\,:\, \PP^2 \times \PP^2 \times \PP^2 \,\longrightarrow \,{\rm Gr}(3,6) \subset \PP^{19}$$
Algebraically, the map $\Phi$ is given by evaluating the twenty
$3 \times 3$ minors of the matrix (\ref{3by6matrix}).

This note concerns the {\em tropicalization} of the map $\Phi$.
We study the following inclusions
\begin{equation}
\label{inclusions}
 {\rm image} (\rm trop(\Phi))
\subset {\rm trop}({\rm image}(\Phi))
\subset
{\rm trop}({\rm Gr}(3,6)) \subset \TP^{19}.
\end{equation}
In general, naive tropicalization does not commute with morphisms; accordingly, we will see that the inclusions in (\ref{inclusions}) are both strict. We know from \cite[\S 5]{SS} and  \cite[Table 2]{HJJS}
that
the {\em tropical Grassmannian} ${\rm trop}({\rm Gr}(3,6))$
has a coarsest fan structure, which is represented as a $3$-dimensional
polyhedral complex with $1005$ maximal polytopes,
namely $990$ tetrahedra and $15$ bipyramids,
which is homologically a bouquet of $126$ 3-spheres.
With the help of {\tt GFan} \cite{gfan}, we computed
the two nested subcomplexes on the left in (\ref{inclusions})
and we found that they are also pure of dimension $3$.
The main point of this note is to furnish the
combinatorial descriptions of these polyhedral complexes which are summarized in Proposition \ref{prop:I_3} and Theorem \ref{thm:IV_3}.

\begin{figure}[h!]
  \begin{center}
  \vskip -0.6cm
  \includegraphics[height=6.2in]{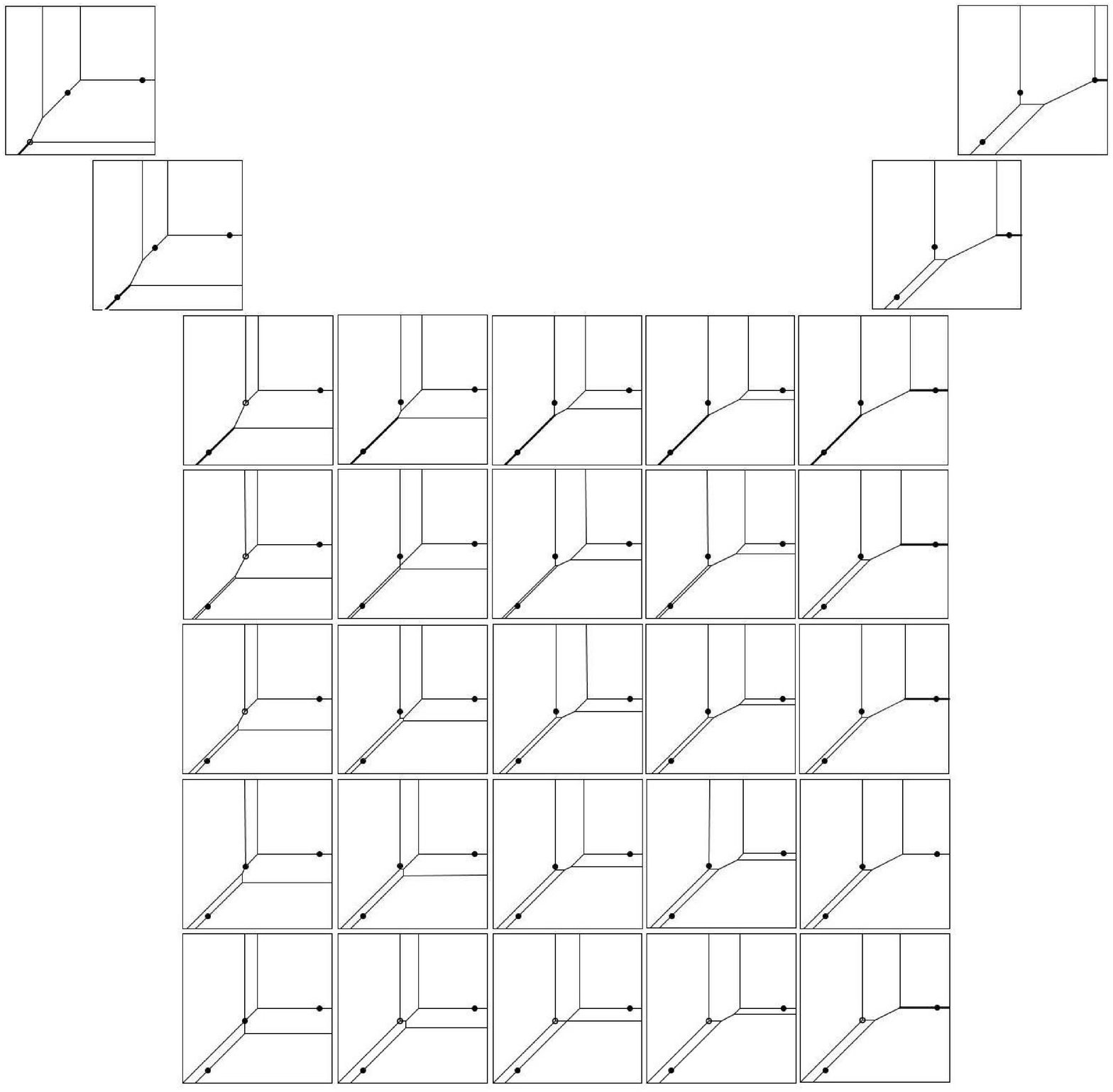}
  \vskip -1cm
  \caption{The plane $L_{X,Y,Z}$ of conics passing through three points
  $X, Y,Z \in \TP^2$.}
    \label{fig:plane1}
  \end{center}
\end{figure}

Many studies in tropical geometry \cite{Mik}
concern curves passing
through given points in the plane $\TP^2$.
Our results complement these by
offering a precise analysis of the plane of conics
passing through three points, as in Figure~\ref{fig:plane1},
 and how that plane depends on the points.

Our results on (\ref{inclusions}) will be stated and derived in Section 3.
In Section 2, we warm up by solving the same problem for
two points in $\PP^2$, where we obtain the {\em tropical Hilbert scheme}
discussed in \cite[\S 6.2]{AN}. Note that
the case of four points in $\PP^2$
was already treated in \cite[\S 6]{RST}.\\

\section{Two Points}

Given two distinct points $x = (x_0{:}x_1{:}x_2)$ and
$y = (y_0{:}y_1{:}y_2)$ in the projective plane $\PP^2$ over a field $K$,
there is a four-dimensional space $L_{x,y}$
of quadrics that vanish at $x$ and $y$, namely
\begin{equation}
\label{2by6matrix}
L_{x,y} \quad = \quad {\rm kernel}
\begin{pmatrix}
x_0^2 & x_1^2 & x_2^2 & x_0 x_1 & x_0 x_2 & x_1 x_2 \\
y_0^2 & y_1^2 & y_2^2 & y_0 y_1 & y_0 y_2 & y_1 y_2
\end{pmatrix}.
\end{equation}
This defines a rational map which takes pairs of points into the Grassmannian:
$$
\Psi \,:\, \PP^2 \times \PP^2 \,\longrightarrow \,{\rm Gr}(4,6) \subset \PP^{14}
$$
The map $\Psi$ blows up the diagonal in
$\PP^2 {\times} \PP^2$. The closure of its image
is isomorphic to the Hilbert scheme
of two points in $\PP^2$. That is,
 $\,{\rm image}(\Psi) = {\rm Hilb}_2(\PP^2) $
    is a smooth $4$-dimensional subvariety of the
  $8$-dimensional   Grassmannian ${\rm Gr}(4,6)$.
  Representing points in ${\rm Gr}(4,6)$ by their dual
  Pl\"ucker coordinates, the map
$\Psi$ is given  algebraically by evaluating the fifteen
$2 {\times} 2$-minors of the matrix in (\ref{2by6matrix}).
Gr\"obner-based implicitization of $\Psi$ yields the prime ideal
$$
\begin{matrix}
I_2  = \! & \!\!\! \!
\langle
p_{03} p_{15} {+} p_{13} p_{34} ,
p_{03} p_{24} {+} p_{04} p_{45} ,
p_{03} p_{25} {+} p_{34} p_{45} ,
p_{04} p_{13} {+} p_{03} p_{35} ,
p_{04} p_{15} {-} p_{34} p_{35} ,
p_{04} p_{25} {-} p_{24} p_{34} , \\ &
p_{13} p_{24} - p_{35} p_{45} ,
p_{13} p_{25} - p_{15} p_{45} ,
p_{15} p_{24} - p_{25} p_{35} ,\,
p_{01} p_{02} - p_{05}^2 + p_{34}^2 ,
p_{01} p_{12} + p_{14}^2 - p_{35}^2 ,\\ &
p_{02} p_{12} - p_{23}^2 + p_{45}^2 , \,
p_{03} p_{23} - p_{34}^2 + p_{04} p_{35} ,
p_{04} p_{14} - p_{34}^2 + p_{03} p_{45} ,
p_{05} p_{15} - p_{35}^2 - p_{13} p_{45} , \\ &
p_{05} p_{25} - p_{24} p_{35} - p_{45}^2 ,\,
p_{13} p_{23} + p_{14} p_{35} - p_{35}^2 - p_{13} p_{45} ,\,
p_{14} p_{24} - p_{24} p_{35} + p_{23} p_{45} - p_{45}^2 ,  \\ & \!\!\!\!\!
p_{01} p_{04} {-} p_{03} p_{05} {-} p_{03} p_{34} ,
p_{01} p_{23} {-} p_{14} p_{34} {+} p_{05} p_{35} ,
p_{01} p_{24} {+} p_{05} p_{45} {+} p_{34} p_{45} ,
p_{01} p_{25} {+} p_{14} p_{45} {+} p_{35} p_{45} , \\ &  \!\!\!\!\!
p_{02} p_{03} {-} p_{04} p_{05} {+} p_{04} p_{34} ,
p_{02} p_{13} {+} p_{05} p_{35} {-} p_{34} p_{35} ,
p_{02} p_{14} {+} p_{23} p_{34} {+} p_{05} p_{45} ,
p_{02} p_{15} {+} p_{23} p_{35} {+} p_{35} p_{45} , \\ & \!\!\!\!\!
p_{03} p_{12} {+} p_{14} p_{34} {-} p_{34} p_{35} ,
p_{03} p_{14} {-} p_{01} p_{34} {+} p_{03} p_{35} ,
p_{04} p_{12} {+} p_{23} p_{34} {-} p_{34} p_{45} ,
p_{04} p_{23} {+} p_{02} p_{34} {+} p_{04} p_{45} ,\\ &  \!\!\!\!\!
p_{05} p_{12} {+} p_{23} p_{35} {-} p_{14} p_{45} ,
p_{05} p_{13} {+} p_{13} p_{34} {+} p_{01} p_{35} ,
p_{05} p_{14} {-} p_{34} p_{35} {+} p_{01} p_{45} ,
p_{05} p_{23} {+} p_{02} p_{35} {+} p_{34} p_{45} ,\\ & \!\!\!\!\!
p_{05} p_{24} {-} p_{24} p_{34} {+} p_{02} p_{45} ,
p_{05} p_{34} {-} p_{04} p_{35} {+} p_{03} p_{45} ,
p_{12} p_{13} {-} p_{14} p_{15} {+} p_{15} p_{35} ,
p_{12} p_{24} {+} p_{23} p_{25} {-} p_{25} p_{45} ,\\ & \!\!\!\!\!
p_{13} p_{14} {+} p_{01} p_{15} {+} p_{13} p_{35} ,
p_{14} p_{23} {+} p_{12} p_{34} {-} p_{35} p_{45} ,
p_{14} p_{25} {-} p_{25} p_{35} {-} p_{12} p_{45} ,
p_{15} p_{23} {+} p_{12} p_{35} {-} p_{15} p_{45} , \\ & \!\!\!\!\!
p_{15} p_{34} {-} p_{14} p_{35} {+} p_{13} p_{45} \,,\,\,
p_{23} p_{24} {+} p_{02} p_{25} {+} p_{24} p_{45} \,,\,\,
p_{25} p_{34} {-} p_{24} p_{35} {+} p_{23} p_{45}
\,\rangle.
\end{matrix}
$$
This is the ideal of the embedding of ${\rm Hilb}_2(\PP^2) $ via $\Psi$
as a subscheme of degree $21$ in  $\PP^{14}$.

We identify the tropicalization of ${\rm Gr}(4,6)$ with the
space of tree metrics on six taxa \cite[\S 2.4]{ascb}.
The taxa are the quadratic monomials.
Combinatorially, this is a simplicial complex
with $25$ vertices, $105$ edges and $105$ triangles.
The $25$ vertices are the splits: trees with  one internal edge.
The tropicalization of $\Psi$ is a piecewise-linear map
into that tree space:
$$
 {\rm trop}(\Psi) \,:\, \TP^2 \times \TP^2 \,\longrightarrow \,{\rm trop}({\rm Gr}(4,6))
\,\subset \, \TP^{14} .
$$
The coordinates of this map are the $15$ tropical $2 {\times} 2$-minors
of the $2 {\times} 6$-matrix in (\ref{2by6matrix}):
$$ {\rm trop}(\Psi)(X,Y) \,=\,
\bigl( {\rm max}(2 X_0 {+} 2Y_1, 2 X_1 {+} 2X_0), \ldots,
{\rm max}(X_0{+}X_2{+}Y_1{+}Y_2, X_1{+}X_2{+}Y_0{+}Y_2) \bigr). $$
We regard its image as a ``combinatorial Hilbert scheme''
that parametrizes pairs of points in $\TP^2$
by the quadrics that pass through them.
We have the following strict inclusions:
\begin{equation}
\label{hilb2inclusions}
\begin{matrix} &
{\rm Hilb}_2(\TP^2) & := & {\rm image}({\rm trop}(\Psi))  & & & \\
\,\subset \, &
{\rm trop}({\rm Hilb}_2(\PP^2)) & =  & {\rm trop}({\rm image}(\Psi))  &
\,\subset \,& {\rm trop} ({\rm Gr}(4,6)) &\,\subset \,& \TP^{14}.
\end{matrix}
\end{equation}
Working modulo the common lineality spaces,
the Hilbert schemes in (\ref{hilb2inclusions}) are one-dimensional complexes.
These graphs are geometrically embedded, but not as subcomplexes,
inside the two-dimensional simplicial complex
 $\,{\rm trop} ({\rm Gr}(4,6))$.
Alessandrini and Nesci argued
in \cite[\S 6.2]{AN} that
${\rm Hilb}_2(\TP^2) $ is a cycle of length six.
The following proposition extends their findings.

\begin{center}
\begin{figure}
\vskip -0.4cm
\includegraphics[width=13.5cm]{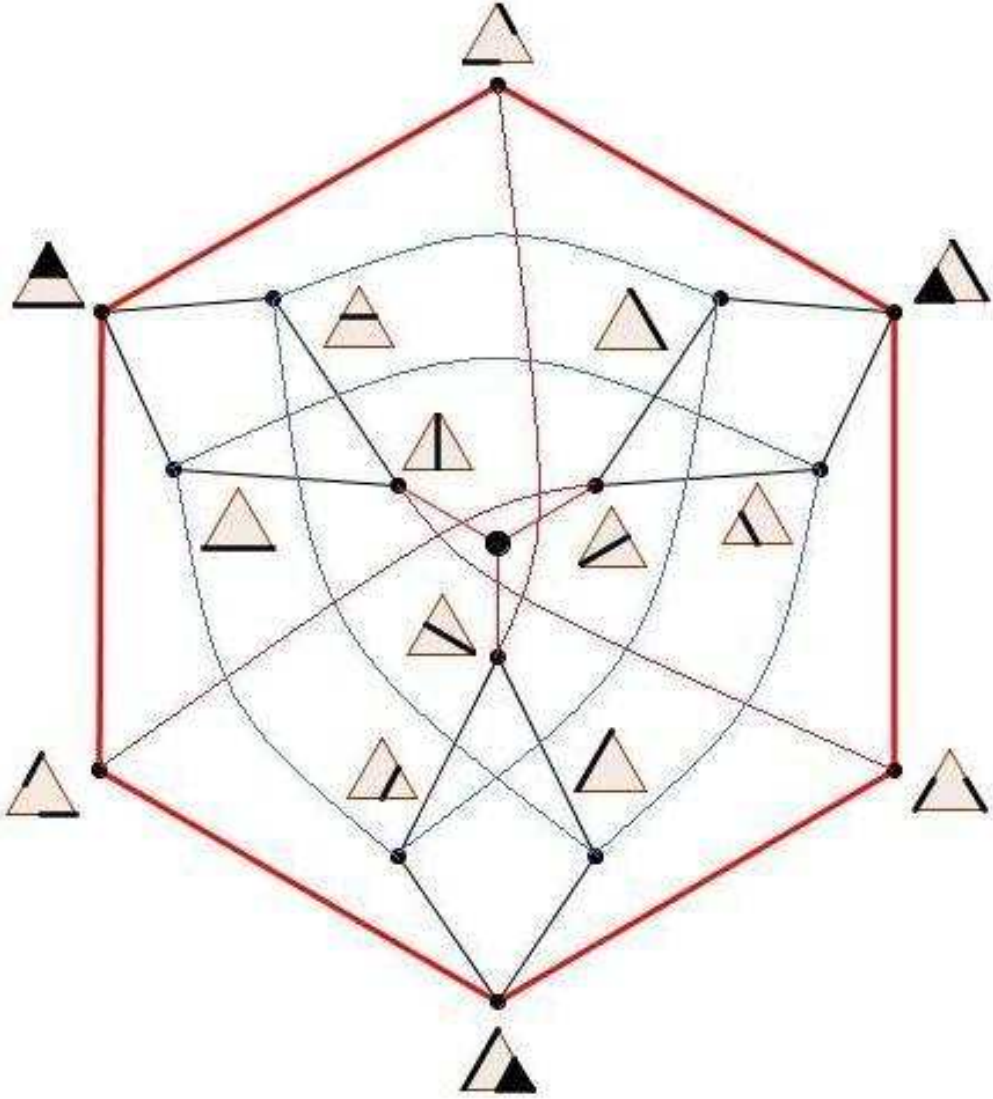}
\vskip -0.2cm
\caption{Tropicalization of the Hilbert scheme of two points in $\PP^2$}
\label{fig:hilb22}
\end{figure}
\end{center}

\begin{proposition}
\label{prop:two}
The tropicalized Hilbert scheme ${\rm trop}({\rm Hilb}_2(\PP^2))$
is the graph with $16$ nodes and $30$ edges  depicted in Figure
\ref{fig:hilb22}. The outer $6$-cycle is the subgraph ${\rm Hilb}_2(\TP^2)$.
The labeling of the graph describes its embedding
in the space of trees on six taxa and is explained below.
\end{proposition}

We proved Proposition  \ref{prop:two} by applying the software {\tt GFan} \cite{gfan}
to the ideal $I_2$ and by carefully analyzing the
output of that computation. We now discuss the
outcome of that analysis.

The graph ${\rm trop}({\rm Hilb}_2(\PP^2))$ has $16$ nodes.
Twelve of the nodes are also nodes in the space of trees,
${\rm trop}({\rm Gr}(2,6))$, so
they correspond to splits of the set
of taxa $\{x_0^2, x_1^2,x_2^2, x_0 x_1, x_0 x_2, x_1 x_2\}$.
Up to the action of the symmetric group $\mathfrak{S}_3$ by permuting the coordinates
of $\PP^2$, there are
\begin{enumerate}
\item  three splits like $\,\{\, \{x_0^2, x_1^2\}\, , \, \{x_2^2, x_0 x_1, x_0 x_2, x_1 x_2\} \,\}$,
\item  three splits like $\,\{\, \{x_0^2, x_1 x_2\}\, , \, \{x_1^2 ,x_2^2, x_0 x_1, x_0 x_2 \} \,\}$,
\item  three splits like $\,\{\, \{x_0 x_1, x_0 x_2\}\, , \, \{x_0^2, x_1^2,x_2^2, x_1 x_2\} \,\}$,
\item  three splits like $\,\{\, \{x_0^2 ,x_0 x_1, x_1^2\}\, , \, \{x_2^2, x_0 x_2, x_1 x_2\} \,\}$.
\end{enumerate}
The nine trees with one split in (1)-(3) have one cherry, or set of edges paired together. We represent
that tree by drawing the cherry pair as a thick black segment in the corresponding
vertex label of Figure~\ref{fig:hilb22}.
The three trees in (4) are 3-3 splits so they have no cherry.
They appear alternatingly on the outer $6$-cycle in Figure~\ref{fig:hilb22},
where they are drawn by a long segment across a triangle.
The other three nodes on the outer $6$-cycle are trees with two interior nodes:
\begin{enumerate}
\item[(5)] three split pairs like
$\bigl\{\! \{\!\{x_0^2,x_0 x_2\}\!, \!\{ x_1^2,x_2^2,  x_0 x_1, x_1 x_2\}\!\},
\{\!\{ x_1^2, x_1 x_2 \} \!, \!\{ x_0^2,x_2^2,  x_0 x_1,  x_0 x_2\}\!\} \! \bigr\}$.
\end{enumerate}
Finally, there is one special node that lies in the relative interior
of a triangle in ${\rm trop}({\rm Gr}(2,6))$:
\begin{enumerate}
\item[(6)] the unique trivalent {\em snowflake tree} with
cherries
$\{ x_0^2, x_1 x_2\}$,
$\{x_1^2, x_0 x_2\}$ and
$\{x_2^2, x_0 x_1\}$.
\end{enumerate}

The graph ${\rm trop}({\rm Hilb}_2(\PP^2))$ has $30$ edges.
Twelve are interior to
 triangles of  ${\rm trop}({\rm Gr}(2,6))$,
 so they correspond to trivalent trees.
Six of those are the edges of type (4-5) that form the
outer $6$-cycle. The others are the three (2-6) edges adjacent to the
snowflake tree, and the three (2-5) edges that appear as the longest edges in Figure \ref{fig:hilb22}.
The remaining $18$ edges of  ${\rm trop}({\rm Hilb}_2(\PP^2))$ are also edges
of ${\rm trop}({\rm Gr}(2,6))$, so they correspond to trees with two interior edges.
Those edges are three (1-2)s, six (1-3)s,
three (1-4)s, three (2-3)s, and three (3-4)s.

The tropical map ${\rm trop}(\Psi)$ amounts to a
double cover of the $6$-cycle ${\rm Hilb}_2(\TP^2)$.
To see this, we note that  the Newton polytope
${\rm NP}(\Psi)$ of the rational map $\Psi$,
as defined in \cite[(3.38)]{ascb} is a centrally-symmetric $12$-gon.
Namely, ${\rm NP}(\Psi)$ is the Minkowski sum of the $15$ Newton polytopes of all
$2 {\times} 2$-minors of the $2 {\times} 6$-matrix (\ref{2by6matrix}), and these are
line segments
that lie in a common plane and involve six distinct directions.
According to \cite[Theorem 3.42]{ascb}, ${\rm trop}(\Psi)$
is linear on each of the twelve cones in the normal fan of ${\rm NP}(\Psi)$.
The $12$-gon formed by these cones is mapped
onto ${\rm Hilb}_2(\TP^2)$ by looping twice around  the outer $6$-gon in Figure~\ref{fig:hilb22}.

\section{Three Points}

Guided by the above results for the
two-point map $\Psi$, we now investigate the
three-point map $\Phi$. We write $I_3$ for the
homogeneous prime ideal that represents the image of
$\Phi$. The following proposition summarizes
basic facts about the variety
$V(I_3) = {\rm image}(\Phi) \subset {\rm Gr}(3,6)$.

\begin{proposition} \label{prop:I_3}
The projective variety $V(I_3)$ has dimension $6$ and degree $57$.
Its ideal  $I_3$ is minimally generated by $62$
homogeneous quadrics in the $20$ Pl\"ucker coordinates $p_{ijk}$.
Among these are $35$ quadrics that vanish on the Grassmannian
${\rm Gr}(3,6) \subset \PP^{19}$.
The corresponding tropical variety ${\rm trop}(V(I_3))$ in $\TP^{19}$,
with its Gr\"obner fan structure
and taken modulo the lineality space, is a $3$-dimensional polyhedral
complex with f-vector $(1095, 6621, 12830, 7649)$.
\end{proposition}

\begin{ex} The quadrics
$\,p_{045} p_{145}-p_{013} p_{245}-p_{345}^2\,$ and
$\, p_{025} p_{145}-p_{134}p_{245}-p_{015} p_{245} + p_{235} p_{345} \,$
are among the $27$ generators of the image of $I_3$
in the coordinate ring of ${\rm Gr}(3,6)$.
\qed
\end{ex}

 There is a natural morphism from the Hilbert scheme
${\rm Hilb}_3(\PP^2)$ onto our variety $V(I_3)$. However, unlike in Section 2,
this is not an isomorphism. Geometrically, the morphism  from the Hilbert scheme
 contracts all triples of points that lie on the same line.
The singular locus of $V(I_3)$ is a projective plane $\PP^2$,
namely, it is the image of all collinear triples in ${\rm Hilb}_3(\PP^2)$.

The Gr\"obner fan structure on the tropical variety in Propsition \ref{prop:I_3}
 is not simplicial: among the $7649$ three-dimensional
polytopes, $876$ have five vertices
($105$ bipyramids and $773$ Egyptian pyraminds),
$27$ have six vertices (all triangular prisms), and $12$ have seven vertices
(two types). The number of
$\mathfrak{S}_3$-orbits of facets of
${\rm trop}({\rm image}(\Phi)) = {\rm trop}(V(I_3))$  is $1318$.

Our results are summarized in Theorem \ref{thm:IV_3}.
The role of the graph in Figure \ref{fig:hilb22}
is now played by the $3$-dimensional complex with $7649$ facets.
It is obviously too big to be fully displayed here. Instead, we shall now focus
on the leftmost complex in (\ref{inclusions}). The tropical morphism ${\rm trop}(\Phi)$ is a piecewise linear map.
As shown in  \cite[\S 3.4]{ascb}, its domains of linearity are the normal cones
of the Newton polytope.
We begin by computing this Newton polytope.

\begin{lemma}
The Newton polytope of $\Phi$ is $4$-dimensional and
has f-vector $(504,1056,684,132)$.
\end{lemma}

The polytope ${\rm NP}(\Phi)$ plays the same role as
the $12$-gon ${\rm NP}(\Psi)$ in Section~2. The $504$ vertices of ${\rm NP}(\Phi)$ correspond
to distinct types of three labeled points in $\TP^2$, where the type is the cell
of ${\rm Trop}({\rm Gr}(3,6))$ that contains the plane of quadrics through these points.
The map $\Phi$ is invariant under permuting the three points $x$, $y$ and $z$, and this reduces
the number of image cones to $504/6 =  84$. These $84$ cones in $\TP^{19}$
are grouped into $17$ orbits of size six with respect to the common symmetry
group $\mathfrak{S}_3$ of $I_3$, $V(I_3)$ and  ${\rm Trop}(V(I_3))$.
Those symmetries correspond to permuting the indices  $0$, $1$ and $2$ of the coordinates on
$\PP^2$ or $\TP^2$.

\begin{table}[h!]
\caption{The $17$ generic types of triples in $\TP^2$ and their planes of conics}.
\centering
\begin{tabular}[h!]{l c c c c r}
\hline
Type & Orbit Size & Valency &  $(X_1,X_2) $ & $ (Y_1,Y_2)$ & Plane $L_{X,Y,Z}$ \\
\hline
1 & 12 & 6 & $(4,3)$ & $(3,-1)$ & FFFGG\\
2  & 36 &  5 & $(6,5)$ & $(3,-1)$  & EFFG\\
3 & 36 & 5 & $(4,5)$ & $(2,-1)$   & FFFGG\\
4 & 36 & 4 &  $(6,7)$ & $(2,-1)$  & EFFG\\
5 & 36 & 4 &  $(8,6)$ & $(5,1)$ & EFFG\\
6 & 36 & 4 &  $(5,8)$ & $(-1,5)$  & EEFG\\
7 & 36 & 4 & $ (5,6)$ & $(2,1)$ & EEFF(a)\\
8 & 36 & 4 &  $(3,7)$ & $(-1,2)$  & EEFG\\
9 & 36 &  4 & $(6,5)$ & $(4,2)$ & EEFF(a)\\
10 & 36 & 4 &  $(5,6)$ & $(3,1)$ & EFFG\\
11 & 36 &  4 & $ (6,5)$ & $(5,2)$ & EFFG\\
12 & 36 &  4 & $(4,6 )$ & $(2,3)$ & EEFF(a)\\
13 & 12 & 4 & $(5,3)$ & $(3,-2)$  & EEEG\\
14 & 18 &  4 & $(3,6)$ & $(1,3)$ & EEFF(b)\\
15 & 18 &  4 &  $(3,6)$ & $(2,3)$ & EEFF(b)\\
16 & 36 &  4 & $(3,5)$ & $(2,1)$ & EEFG\\
17 & 12 &  4 & $(2,3)$ & $(3,1)$ & EEEG\\
\hline
\end{tabular}
\label{table:maxcones}
\end{table}

Given three points $X = (X_0,X_1,X_2)$,
$Y = (Y_0,Y_1,Y_2)$ and $Z = (Z_0,Z_1,Z_2)$ in the tropical projective plane $\TP^2$, we
write $L_{X,Y,Z}$ for the tropical $2$-plane in $\TP^5$ determined,
as in \cite[(3.44)]{ascb} or
 \cite[\S 2]{Sp}, by the tropical Pl\"ucker vector
$\,{\rm trop}(\Phi)(X,Y,Z) \,\in \,{\rm trop}({\rm Gr}(3,6)) \,\subset \,\TP^{19}$.
Geometrically, $L_{X,Y,Z}$ is the tropical plane whose points are the
tropical quadrics that pass through the points $X,Y,Z$.
Our picture of this in Figure~\ref{fig:plane1} is
 reminiscent of \cite[Fig.~19]{RST}.
 
 Our main result is the classification of the $17$ types
of configurations of triples of points.

\begin{theorem}\label{thm:IV_3}
Precisely $48$ of the $1005$ generic $2$-planes in $\TP^5$ arise
as $L_{X,Y,Z}$ for some triple $X,Y,Z \in \TP^2$.
This covers six of the seven symmetry classes.
Table~\ref{table:maxcones} summarizes
the correspondence between the $17$ types of triples and
the $6$ combinatorial types of $2$-planes.
\end{theorem}

We proved this theorem by explicit computations. We shall explain our method and how to read Table~\ref{table:maxcones}.
For each of the $17$, types we list a representative configuration.
Here we break the symmetry by setting $X_0 = Y_0 = Z_0 = 0$
and by fixing the third point to lie at the origin, i.e.~$Z = (Z_1,Z_2) = (0,0)$.
For each configuration, the first point $X= (X_1,X_2)$ is listed in the fourth column,
and the second point $Y=(Y_1,Y_2)$ is listed in the fifth column.

\begin{figure}
%\begin{figure}[h!]
  \begin{center}
  \includegraphics[height=4.3in]{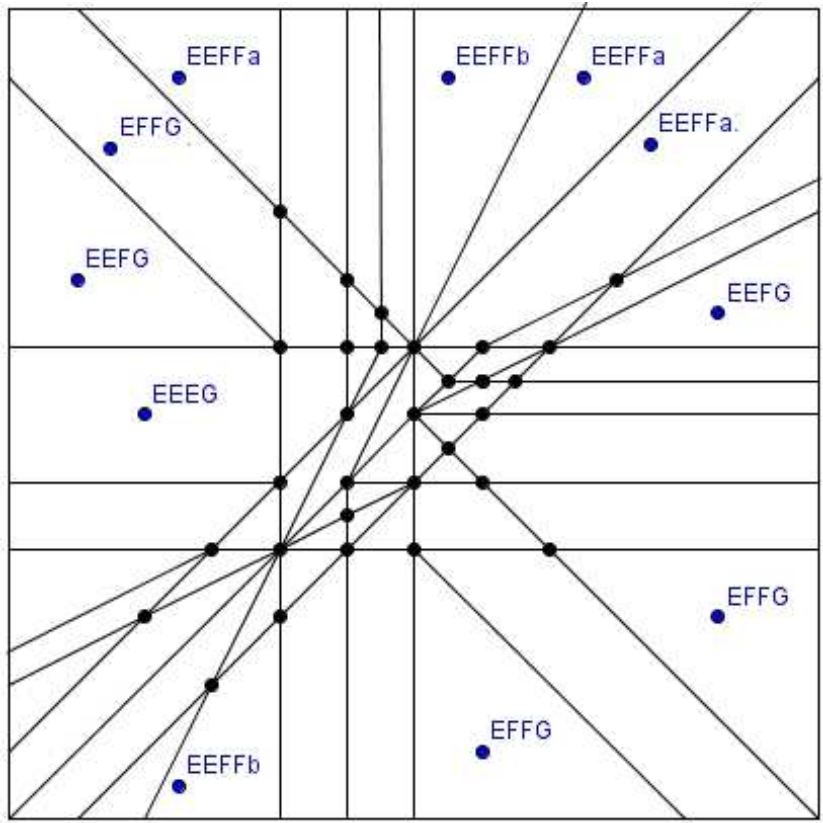}\\
  \caption{Partition of the $(Z_1,Z_2)$-plane obtained
  by $X = (0,0)$ and $Y= (2,3)$.}
    \label{DivisionPlane}
  \end{center}
  % Requires \usepackage
\end{figure}

The second column, ``Orbit Size'', lists of the cardinality of the orbit of the configuration
under permuting both points and coordinates. The sum of the $17$ orbit
sizes is $504$, the total number of vertices of ${\rm NP}(\Phi)$. The third column, ``Valency'', lists the number of edges of
the Newton polytope ${\rm NP}(\Phi)$ that are adjacent to the
given vertex. Equivalently, this is the number of linear inequalities needed to
characterize the configurations of the type in question.
For instance, there are
six such linear inequalities for type 1:
\begin{equation}
\label{eq:sixineq}
X_1 \geq X_2, \,\, X_1 \geq Y_1,\,\,
 Y_2 \leq 0 ,\,\,
  X_2+2 Y_1 \geq 2 X_1, \,\,
    2 X_2+Y_1 \geq 2 X_1,\,\,
      X_2+Y_1+Y_2 \geq X_1.
 \end{equation}
 One solution is $\{X = (4,3), Y = (3,-1) \}$. A solution is equivalent,
  in the sense that the plane
$L_{X,Y,Z}$ is in the same maximal cone
of ${\rm trop}({\rm Gr}(3,6))$, if and only if the six inequalities in
(\ref{eq:sixineq}) are satisfied.
The solution set of (\ref{eq:sixineq}) is the
cone over a bipyramid. The corresponding vertex figure of
${\rm NP}(\Phi)$ is a cube. The vertex figures for
types 2 and 3 are Egyptian pyramids.
All other $14$ types of vertices are simple, so those vertex
figures of ${\rm NP}(\Phi)$ are tetrahedra.

One way to visualize the partition of $\R^4$ into $504$ normal
cones to ${\rm NP}(\Phi)$ is to intersect this normal fan with the
$2$-dimensional affine space obtained by also fixing
the second point $Y = (Y_1,Y_2)$
at a particular location. For instance, in Figure \ref{DivisionPlane}
we fix  $Y=(2,3)$, and we allow $X = (X_1,X_2)$ to vary over
the plane. The regions of equivalence are convex polygons.
On each polygon, the plane of conics $L_{X,Y,Z}$ has a fixed
combinatorial type in ${\rm trop}({\rm Gr}(3,6))$.\

The tropicalization of ${\rm Gr}(3,6)$ has seven $\mathfrak{S}_6$-orbits of maximal cones. The interior of each corresponds to a distinct type
of plane in $\TP^5$. These types were classified in \cite[\S 5]{SS} and  given the names EEEE, EEFF(a), EEFF(b), EFFG, EEEG, EEFG, and FFFGG.
See \cite[Figure 1]{HJJS} for a diagram that shows these seven planes.
In the last column of Table~\ref{table:maxcones}, we see that
 EEEE is the unique type that does not arise as $L_{X,Y,Z}$
for any triple $X,Y,Z$ in $\TP^2$. The other six types arise
as planes of conics through three points, as also seen
in Figure \ref{DivisionPlane}.

Each tropical plane consists of bounded and unbounded faces.
These are dual to the interior cells and boundary cells of
a matroid subdivision of the second hypersimplex \cite{Sp}.
Each cell is indexed by a matroid of rank $3$ on
the ground set $\{1,2,3,4,5,6\}$. Intersecting a generic tropical plane
in $\TP^5$ with the six tropical hyperplanes at infinity gives a tree arrangement consisting of six trees with five leaves each. A detailed description of the correspondence between tropical planes
and tree arrangements is given in \cite[\S 4]{HJJS}. For our planes $L_{X,Y,Z}$ that arise from triples $X,Y,Z \in \TP^2$,
we can construct the corresponding arrangement of six trees by removing in turn each of the
six monomials $\{U_0^2, U_1^2, U_2^2, U_0 U_1, U_0 U_2, U_1 U_2 \}$
from the defining equation of the conic. In other words, each of the
six trees represents a line in $\TP^4$. That line is a tree which parameterizes
conics with five fixed terms that pass through the three given points.
These trees are similar to  \cite[Fig.~19]{RST} but  have only five taxa.\

Now we explain how we constructed Table \ref{table:maxcones}.
We first computed the Newton polytope ${\rm NP}(\Phi)$
using {\tt Gfan}, and we picked a representative in each
maximal cone of its normal fan. Up to symmetries, these
are the $17$ displayed configurations $X = (X_1,X_2), Y= (Y_1,Y_2) , Z = (0,0)$.
For these we calculated the corresponding vectors
of tropical Pl\"ucker coordinates
$$ ({\rm trop}(\Phi))(X,Y,Z) \,\,= \,\,
\bigl( P_{012}, P_{013}, P_{014}, \ldots, P_{245}, P_{345} \bigr) \,\, \in \,\, \R^{20}. $$
For each of the six indices, we consider the restricted vector
of Pl\"ucker coordinates involving that index. For instance, for index
``$0$'', corresponding to the monomial $U_0^2$, this is the vector
$$(P_{012}, P_{013}, P_{014}, P_{023}, P_{024}, P_{034}) .$$
This vector represents the pairwise distances in a phylogenetic
tree with taxa $\{1,2,3,4,5\}$. This is the first among the six trees
that represent the plane $L_{X,Y,Z}$, in its guise as
a $3$-tree \cite[(3.44)]{ascb}. At this point, the last column in
Table \ref{table:maxcones} can simply read off from
\cite[Table 2]{HJJS}.

\bigskip

{\bf Acknowledgements.}
We are grateful to  Anders Jensen for helping us with our {\tt Gfan} computation.
This work is based on the undergraduate senior thesis of the first author.
The second author was partially supported by the NSF
 (DMS-0456960 and DMS-0757207).


\bigskip


\begin{thebibliography}{99}

\bibitem{AN} D.~Alessandrini and M.~Nesci:
On the tropicalization of the Hilbert scheme,
{\tt arXiv:0912.0082}.

\bibitem{HJJS} S. Herrmann, A. Jensen, M. Joswig, and B. Sturmfels:
How to draw tropical planes, {\em Electron.~J.~Combin.} {\bf 16} (2009) \# 6.

\bibitem{gfan} A Jensen: Gfan, a Software System for {G}r{\"o}bner fans and tropical varieties, Available at http://www.math.tu-berlin.de/\~{}jensen/software/gfan/gfan.html.

\bibitem{Mik} G.~Mikhalkin: Tropical geometry and its applications. International Congress of Mathematicians. Vol. II, 827--852, Eur. Math. Soc., Z\"urich, 2006.

\bibitem{ascb} L.~Pachter and B.~Sturmfels:
{\em Algebraic Statistics for Computational Biology},
Cambridge University Press, 2005.

\bibitem{RST}
J.~Richter-Gebert, B.~Sturmfels and T.~Theobald:
First steps in tropical geometry.
Idempotent mathematics and mathematical physics, 289--317,
{\em Contemp. Math.}, {\bf 377}, Amer.Math.Soc., Providence, 2005.

\bibitem{Sp} D.~Speyer:
Tropical linear spaces, {\em SIAM J. Discrete Math.} {\bf 22} (2008) 1527--1558.

\bibitem{SS} D.~Speyer and B.~Sturmfels:
The tropical Grassmannian, {\em Advances in Geometry} {\bf 4} (2004) 389--411.

\end{thebibliography}
\end{document}